\documentclass[12pt]{amsart}
\usepackage{amsmath}

\newcommand{\nc}{\newcommand}
\nc{\cal}{\mathcal} 
\nc{\la}{\langle}
\nc{\ra}{\rangle}
 \nc{\CA}{\cal A}
 \nc{\CBB}{\cal B}
 \nc{\CC}{\cal C}
\nc{\CDD}{\cal D}
\nc{\CE}{\cal E}
\nc{\CF}{\cal F}
\nc{\CG}{\cal G}
\nc{\CH}{\cal H}
\nc{\CI}{\cal I}
\nc{\CJ}{\cal J}
\nc{\CK}{\cal K}
\nc{\CL}{\cal L}
\nc{\CM}{\cal M}
\nc{\CN}{\cal N}
\nc{\CO}{\cal O}
\nc{\CP}{\cal P}
\nc{\CQ}{\cal Q}
\nc{\CR}{\cal R}
\nc{\CS}{\cal S}
\nc{\CT}{\cal T}
\nc{\CU}{\cal U}
\nc{\CV}{\cal V}
\nc{\CW}{\cal W}
\nc{\CZ}{\cal Z}

\nc{\fa}{\mathfrak a}
\nc{\fg}{\mathfrak g}
\nc{\fk}{\mathfrak k}
\nc{\fh}{\mathfrak h}
\nc{\fm}{\mathfrak m}
\nc{\fn}{\mathfrak n}
\nc{\fA}{\mathfrak A}
\nc{\fC}{\mathfrak C}
\nc{\fI}{\mathfrak I}
\nc{\fL}{\mathfrak L}
\nc{\fS}{\mathfrak S}

\nc{\nen}{\newenvironment}
\nc{\ol}{\overline}
\nc{\ul}{\underline}
\nc{\lra}{\longrightarrow}
\nc{\lla}{\longleftarrow}
\nc{\Lra}{\Longrightarrow}
\nc{\Lla}{\Longleftarrow}
\nc{\Llra}{\Longleftrightarrow}
\nc{\hra}{\hookrightarrow}
\nc{\iso}{\overset{\sim}{\lra}}
\nc{\Hom}{\mathrm{Hom}}
\nc{\Mor}{\mathrm{Mor}}

\nc{\notebox}[1]{\noindent\fbox{\parbox{12.5cm}{\sf #1}}\\[8pt]}
\nc{\Thm}[1]{Theorem~\ref{#1}}
\nc{\Prop}[1]{Proposition~\ref{#1}}
\nc{\Lem}[1]{Lemma~\ref{#1}}
\nc{\Cor}[1]{Corollary~\ref{#1}}
\nc{\Conj}[1]{Conjecture~\ref{#1}}
\nc{\Claim}[1]{Claim~\ref{#1}}
\nc{\Defn}[1]{ Definition~\ref{#1}}
\nc{\Exa}[1]{Example~\ref{#1}}
\nc{\Rem}[1]{Remark~\ref{#1}}
\nc{\Note}[1]{Note~\ref{#1}}

\nc{\marg}{\marginpar}

\nen{thm}[1]{\label{#1}{\bf Theorem.\ } \em}{}
\nen{thma}[1]{\label{#1}{\bf Theorem A.\ } \em}{}
\nen{thmb}[1]{\label{#1}{\bf Theorem B.\ } \em}{}
\nen{prop}[1]{\label{#1}{\bf Proposition.\ } \em}{}
\nen{lem}[1]{\label{#1}{\bf Lemma.\ } \em}{}
\nen{klem}[1]{\label{#1}{\bf Key Lemma.\ } \em}{}
\nen{cor}[1]{\label{#1}{\bf Corollary.\ } \em}{}
\nen{cora}[1]{\label{#1}{\bf Corollary A.\ } \em}{}
\nen{coraa}[1]{\label{#1}{\bf Corollary A'.\ } \em}{}
\nen{corb}[1]{\label{#1}{\bf Corollary B.\ } \em}{}
\nen{corbb}[1]{\label{#1}{\bf Corollary B'.\ } \em}{}
\nen{conj}[1]{\label{#1}{\bf Conjecture.\ } \em}{}
\nen{conja}[1]{\label{#1}{\bf Conjecture A.\ } \em}{}
\nen{conjb}[1]{\label{#1}{\bf Conjecture B.\ } \em}{}
\nen{claim}[1]{\label{#1}{\bf Claim.\ } \em}{}

\nen{defn}[1]{\label{#1}{\bf Definition.\ } }{}
\nen{exa}[1]{\label{#1}{\bf Example.\ } }{}

\nen{rem}[1]{\label{#1}{\em Remark.\ } }{}
\nen{exer}[1]{\label{#1}{\em Exercise.\ } }{}
\nen{pf}{\begin{proof}}{\end{proof}}

 \nc{\br}{\mathbb R}
 \nc{\bz}{\mathbb Z}
 \nc{\bc}{\mathbb C}
 \nc{\bn}{\mathbb N}
 \nc{\geg}{\mathfrak g}
 \nc{\G}{\Gamma}
 \nc{\sm}{\setminus}
 \nc{\sub}{\subset}
 \nc{\lm}{\lambda}
 \nc{\al}{\alpha}
  \nc{\bt}{\beta}
 \nc{\om}{\omega}
  \nc{\dl}{\delta}
  \nc{\g}{\gamma}
 \nc{\Dl}{\Delta}
 \nc{\Om}{\Omega}
 \nc{\s}{\sigma}
 \nc{\ro}{\rho}
  \nc{\te}{\theta}
 \nc{\SLR}{SL_2(\br)}
 \nc{\GLR}{GL_2(\br)}
 \nc{\PGLR}{PGL_2(\br)}
 \nc{\PSLR}{PSL_2(\br)}
 \nc{\SLC}{SL(2,\bc)}
 \nc{\uH}{\mathbb H}
 \nc{\fD}{\mathcal{D}}
 \nc{\fE}{\mathcal{E}}
 \nc{\fO}{\mathcal{O}}
 \nc{\haf}{\frac{1}{2}}
 \nc{\qtr}{\frac{1}{4}}
 \nc{\8}{\infty}
  \nc{\7}{{-\infty}}
 \nc{\inv}{^{-1}}
 \nc{\eps}{\varepsilon}
\begin{document}
\title[restrictions of eigenfunctions]
{Estimates of geodesic restrictions of eigenfunctions on
hyperbolic surfaces and representation theory {\footnote{\tiny{ \it This is an updated version of the text not intended for a publication for being obsolete.}}}}

\author{Andre Reznikov}
\address{Bar-Ilan University, Ramat Gan, Israel}
\email{reznikov@math.biu.ac.il}
\begin{abstract} We consider restrictions along closed geodesics
and geodesic circles for eigenfunctions of the Laplace-Beltrami
operator on a compact hyperbolic Riemann surface. We obtain bounds
on the $L^2$-norm and on the generalized periods of such
restrictions as the corresponding eigenvalue tends to infinity. We
use methods from the theory of automorphic functions and in
particular the uniqueness of the corresponding invariant
functionals on irreducible unitary representations of $\PGLR$.
\end{abstract}
\maketitle

\section{Introduction}
\label{intro}

\subsection{Maass forms}\label{M-forms}
Let $Y$ be a compact Riemann surface  with a Riemannian metric of
constant curvature $-1$. We denote by $dv$ the associated volume
element and by ${\rm d}(\cdot,\cdot)$ the corresponding distance
function. The corresponding Laplace-Beltrami operator $\Dl$ is
non-negative and has purely discrete spectrum on the space
$L^2(Y,dv)$ of functions on $Y$. We will denote by $0=\mu_0< \mu_1
\leq \mu_2 \leq ...$ the eigenvalues of $\Dl$ and by $\phi_i$ the
corresponding eigenfunctions (normalized to have $L^2$ norm one).
In the theory of automorphic functions the functions $\phi_i$ are
called non-holomorphic forms, Maass forms (after H. Maass,
\cite{M}) or simply automorphic functions.

The study of Maass forms is important in analysis, number theory
and in many areas of mathematics and of mathematical physics. In
particular, various questions concerning analytic properties of
the eigenfunctions $\phi_i$ drew a lot of attention in recent
years (see surveys \cite{Sa1}, \cite{Sa2} and references therein).

\subsection{Restrictions to curves}\label{periods}
In this paper we study analytic properties of restrictions of
eigenfunctions to closed curves.

Let $\gamma\subset Y$ be a smooth closed curve with the
corresponding line element $d\g$. We fix a parametrization
$t_\g:S^1\to \gamma$ which we assume is proportional to the
natural parametrization (i.e. $d\g^*= length(\g)d\te$, where
 $0\leq\te< 1$ is the parametrization of $S^1$). Our
main object of the study is the restriction
$\phi_i^\g(\te)=\phi_i(t_\g(\te))$ of a Maass form $\phi_i$ to the
curve $\g$.

The first question naturally arising in this setting asks for a
bound on the size of the $L^2$-norm of the restriction of $\phi_i$
to $\g$:
\begin{align}\label{per-def}
p^\g(\phi_i)=\int_\g|\phi_i|^2d\g\
\end{align}
as the eigenvalue $\mu_i\to\8$.

We study the $L^2$-norm of such a restriction via generalized or
twisted periods of a Maass form along the closed parameterized
curve. Such periods are naturally arise in the theory of
automorphic functions and are of interest in their own right.
Namely, as a function on $S^1$, the restriction $\phi_i^\g$ gives
rise to the following Fourier coefficients :
\begin{align}\label{g-per-def}
p_n^\g(\phi_i)=\int_\g\phi_i(t_\g(\te))e^{-2\pi in\te}d\g\ .
\end{align}

We consider two types of curves on the Riemann surface $Y$: closed
geodesics and geodesic circles. It turns out that for these
special curves one can study the corresponding Fourier
coefficients of Maass forms via representation theory.

\subsubsection{Geodesic circles}\label{cir-geo}
We consider geodesic circles first. Let $\s=\s(r,y)\subset Y$ be a
geodesic circle of a radius $r>0$ centered at $y\in Y$. For $r$
which is less than the injectivity radius of $Y$ at $y$ the
geodesic circle is defined by $\s(r,y)=\{y'\in Y|\ {\rm
d}(y,y')=r\}$. Writing the Laplace-Beltrami operator in polar
geodesic coordinates centered at $y$ and using the separation of
variables one sees immediately that there exists a function
$C_\mu(r,n)$ such that the $n$'th generalized period of a Maass
form $\phi_i$ with the eigenvalue $\mu_i$ along $\s(r,y)$ is equal
to
\begin{align}\label{c-period-normal}
p^{\s(r,y)}_n(\phi_i)=\int_{\s(r,y)}\phi_i(t_\s(\te))e^{-2\pi
in\te}d\s=a^\s_n(\phi_i)C_{\mu_i}(r,n) .
\end{align}
In other words, the restriction of the Maass form $\phi_i$ to the
geodesic circle $\s(r,y)$ has the Fourier series expansion given
by
\begin{align}\label{c-fourier}
\phi_i(t_\s(\te))=\sum_n a^\s_n(\phi_i)C_{\mu_i}(r,n)e^{2\pi
in\te}\ , \ n\in \bz .
\end{align}

 We stress that the function $C_\mu(r,n)$ depends only on the
eigenvalue of $\phi_i$ and not on the choice of the eigenfunction
(it is essentially equal to the appropriate hypergeometric
function or Legendre function, see \cite{He}, \cite{Sa3},
\cite{Pe}) and is independent of the point $y$. On the other hand
coefficients $a^\s_n(\phi_i,y)$ capture the structure of the
eigenfunction $\phi_i$. The expansion \eqref{c-fourier} is similar
to the Taylor expansion at the point $y\in Y$, and in fact would be the Taylor expansion for the holomorphic forms.

Coefficients $a^\s_n(\phi_i,y)$ are the main object of our study.
Our main result is the following

\begin{thma}{main-thm2}
For any fixed geodesic circle $\s=\s(r,y)\subset Y$ there exists a
constant $C_\s$ such that for any eigenfunction $\phi_i$ with the
eigenvalue $\mu_i$ the following bound holds
\begin{align}\label{main-ave-bound-K}
\sum _{|n|\leq T}|a^\s_n(\phi_i)|^2\leq C_\s\cdot
\max\{T,\sqrt{\mu_i}\}\ ,
\end{align}
for any $T\geq 1$.
\end{thma}

In fact, we will prove this theorem for $\s$ being any {\it fixed}
image (not necessary smooth) of a circle in a tangent space $T_yY$
at $y$ under the exponential map $\exp_y:T_yY\to Y$.

As a corollary we obtain the following bound on the $L^2$-norm of
the restriction.

\begin{cora}{cor-A} Under the same conditions as in the Theorem A
there exists a constant $C'_\s$ such that the following bound
holds
\begin{align}\label{main-bound-K}
p^\s(\phi_i)\leq C'_\s\cdot \mu_i^{\frac{1}{6}}\ .
\end{align}
\end{cora}
It is natural to expect that the true value of the norm for the restriction is given by the following

\begin{conja}{conjA}
Let $\s\subset Y$ be a fixed closed geodesic circle then
\begin{align}
p^\s(\phi_i)\ll \mu_i^\eps
\end{align}for any $\eps>0$.
\end{conja}

We also obtain a uniform bound for the period.

\begin{coraa}{cor-A2} Under the same conditions as in the Theorem A
there exists a constant $C''_\s$ such that the following bound
holds
\begin{align}\label{period-bound-K}
|p^{\s(r,y)}_0(\phi_i)|=\left|\int_{\s(r,y)}\phi_i(t_\s(\te))d\s\right|\leq C''_\s\ .
\end{align}
\end{coraa}

In fact, one can formulate a similar bound for the generalized period $p^{\s(r,y)}_0(\phi_i)$ for any $n$ (see Section \ref{K-period}). On the basis of the Lindelof conjecture, one expects that the bound $|\int_{\s(r,y)}\phi_i(t_\s(\te))d\s|\ll \mu_i^{-1/4+\eps}$ holds, but no improvement over  \eqref{period-bound-K} is known for general surfaces. For arithmetic surfaces and special circles (coming from imaginary quadratic fields) non-trivial improvements follow from the subconvexity bounds for the corresponding $L$-functions (see \cite{MV}). We note that the analogous to \eqref{period-bound-K} bound is sharp on the sphere $S^2$ and on the torus $T^2$.

\subsubsection{Closed geodesics}\label{cl-geo}
We have a similar statement for closed geodesics as well. Namely,
there exists a function $G_\mu(l,n)$ such that for any closed
geodesic $\ell$ of a length $l$ the corresponding period is given
by
\begin{align}\label{g-period-normal}
p^{\ell}_n(\phi_{\mu_i})=\int_{\ell}\phi_i(t_\ell(\te))e^{-2\pi
in\te}d\ell=a^\ell_n(\phi_i)G_{\mu_i}(l,n) .
\end{align}

The existence of the function $G_\mu(l,n)$ again follows from the
separation of variables in coordinates corresponding to the first
coordinate being the natural parameter along the geodesic $\ell$
and the second coordinate being the distance to the geodesic
$\ell$. It is also can be described in terms of an appropriate
hypergeometric function. In Section \ref{K-restrict} we will show
how the existence of the function $G_\mu(l,n)$ easily follows from
the representation theory.

\begin{thmb}{main-thm}
For any closed geodesic $\ell\subset Y$ there exists a constant
$C_\ell$ such that for any eigenfunction $\phi_i$ with the
eigenvalue $\mu_i$ the following bound holds
\begin{align}\label{main-ave-bound-A}
\sum _{|n|\leq T}|a^\ell_n(\phi_i)|^2\leq C_\ell\cdot
\max\{T,\sqrt{\mu_i}\}\ ,
\end{align}
for any $T\geq 1$.

\end{thmb}
 We have the similar

\begin{corb}{cor-b} Under the same conditions as in the Theorem B
there exists a constant $C'_\ell$ such that the following bound
holds
\begin{align}\label{main-bound}
p^\ell(\phi_i)\leq C'_\ell\cdot \mu_i^{\frac{1}{4}}\ .
\end{align}
\end{corb}

The reason for the difference in exponents in Corollaries A and B
is explained in Sections \ref{K-restrict} and \ref{K-period}.

It is natural to expect that the true value of the norm for the restriction is given by the following

\begin{conjb}{conjB}
Let $\ell\subset Y$ be a fixed closed geodesic then
\begin{align}
p^\ell(\phi_i)\ll \mu_i^\eps
\end{align}for any $\eps>0$.
\end{conjb}

We also obtain a uniform bound for the period (compare to \cite{Pitt}).

\begin{corbb}{cor-B2} Under the same conditions as in the Theorem A
there exists a constant $C''_\ell$ such that the following bound
holds
\begin{align}\label{period-bound-K}
|p^{\ell}_0(\phi_i)|=\left|\int_{\ell}\phi_i(t_\ell(\te))d\s\right|\leq C''_\ell\ .
\end{align}
\end{corbb}

In fact, one can formulate a similar bound for the generalized period $p^{\ell}_0(\phi_i)$ for any $n$ (see Section \ref{K-restrict}). On the basis of the Lindelof conjecture, one expects that the bound $|\int_{\ell}\phi_i(t_\s(\te))d\s|\ll \mu_i^{-1/4+\eps}$ holds, but no improvement over  \eqref{period-bound-K} is known for general surfaces. For arithmetic surfaces and special circles (coming from imaginary quadratic fields) non-trivial improvements follow from the subconvexity bounds for the corresponding $L$-functions (see \cite{MV}). We note that the analogous to \eqref{period-bound-K} bound is sharp on the sphere $S^2$ and on the torus $T^2$.

\subsection{The method}\label{method}
We follow the general strategy formulated in \cite{BR3}. It is
based on ideas from the representation theory of the group
$\PGLR$. Namely, any Riemannian surface $Y$ gives rise to a
discrete subgroup $\G\subset\PGLR=G$ and the quotient space
$X=\G\sm G$. We first use the fact that every Maass form $\phi$
generates an irreducible unitary representation $V\subset L^2(X)$
of $G$, called an automorphic representation. The eigenfunction
$\phi$ corresponds to a unit vector $e_0\in V$ invariant under the
compact subgroup $K\subset G$ such that $X/K\simeq Y$ (such a
vector is unique up to multiplication by a constant). All
irreducible unitary representations of $G$ are classified and have
{\it explicit models}. This setting, formulated by I. Gel'fand and
S. Fomin (\cite{GF}), will be essential in what follows.

Generalized periods along closed geodesics and geodesic circles
are related to the representation theory in the following way. Any
closed geodesic $\ell$ gives rise to a closed orbit $\fO_\ell$ of
the diagonal subgroup $A\subset G$ under the right action of $G$
on $X$ and similarly each geodesic circle $\s$ gives rise to an
orbit $\fO_\s$ of an appropriate compact subgroup $K'$ of $G$. We
note that any closed orbit $\fO$ of $A$ gives rise to the cyclic
subgroup $A_\fO\subset A$ of elements in $A$ acting trivially on
$\fO$. The quotient group $A/A_\fO$ is compact (this is also true
for an orbit of a compact subgroup but we will not use this
subgroup since $K'$ is compact by itself). The $L^2$-form on the
orbit $\fO$ gives rise to a non-negative Hermitian form $H_\fO$ on
the space of smooth functions on $X$ via restriction. We study the
coefficients $p^{\ell/\s}(\phi_i)$ and $p_n^{\ell/\s}(\phi_i)$
through this Hermitian form.

Let $\mathcal{G}\subset G$ stands for either $A$ or $K'$. For an
orbit $\fO$ of $\mathcal{G}$ as above we study the form $H_\fO$ by
means of the corresponding generalized periods along $\fO$.
Namely, for a unitary character $\chi:\mathcal{G}\to\bc$ we
consider a functional on the space $V\cap C^\8(X)$, which we call
the $\chi$-period along the orbit $\fO$, defined on $V\subset
C^\8(X)$ by via the integral against $\bar\chi$ over $\fO$ with
respect to a $\mathcal{G}$-invariant measure. For $\mathcal{G}=A$
we consider only characters trivial on $A_\fO$. One have the
Plancherel formula expressing the value of the form $H_\fO$ on a
vector in $V$ as the sum of squares of the coefficients arising
from $\chi$-periods of the same vector (see \eqref{Plancherel-1}).
This corresponds to the usual Plancherel formula for the Fourier
expansion of the restriction of a function to the orbit $\fO$
under the identification $\fO\simeq S^1$. On the other hand, any
$\chi$-period restricted to an automorphic representation $V$
gives rise to a functional $d^{aut}_\chi$ which is
$\chi$-equivariant under the action of $\mathcal{G}$ on $V$. As
well-known in the representation theory of $\PGLR$  the space of
such functionals on a unitary irreducible representation of $G$ is
one-dimensional. We use this fact in order to introduce the
coefficients $a_\chi\in \bc$ defined below, which are the main
object of our study. Namely, using an {\it explicit model} of an
irreducible automorphic representation in the space of functions
on the line (for $\mathcal{G}=A$) or on the circle (for
$\mathcal{G}=K'$) we define a $(\chi,\mathcal{G})$-equivariant
functional $d_\chi^{mod}$ by means of an {\it explicit kernel}
(see \eqref{coef-a-def}). The uniqueness of such functionals
implies the existence of the coefficient of proportionality
$$d^{aut}_\chi=a_\chi\cdot d^{mod}_\chi\ .$$
The main idea of our method is to use analytic properties of the
explicitly constructed functionals $d_\chi^{mod}$ in order to
control the coefficients in the Plancherel decomposition of
$H_\fO$. The study of the functionals $d_\chi^{mod}$ is based on
the stationary phase method. We use simple geometric properties of
the form $H_\fO$ (see Lemmas \ref{L-2-bound} and
\ref{L-2-bound-K}) in order to obtain sharp on the average bound
on the coefficients $a_\chi$ (Theorems \ref{a-bound-thm} and
\ref{a-bound-thm-K}) and use it to prove Theorems A and B. The
discrepancy in the exponent in Corollaries A and B is a reflection
of different spectral decomposition of the same $K$-fixed vector
$e_0\in V$ with respect to characters of $A$ and of $K'$
respectively (compare \eqref{K-value-mod} to
\eqref{K-spec-bound-K}; see Remark \ref{b-ver-c}).

\subsection{Remarks}

1. The supremum norm of an eigenfunction $\phi_\mu$ of the
Laplace-Beltrami operator on a compact Riemannian manifold $M$,
$\dim M=n$, satisfies the H\"{o}rmander's classical  bound
$\sup|\phi_\mu|\leq c\mu^{\frac{n-1}{4}}\|\phi_\mu\|_{L^2(M)}$,
where $\mu$ is the eigenvalue of $\phi_\mu$. This bound is sharp
on the standard sphere. Hence bounds in Corollaries A and B do not
follow from the general pointwise bound on $\phi_i$.

After announcing the proof of bounds in Corollaries A and B author
have learned that these results are special cases of results of D.
Tataru (\cite{Ta}) and also follow from recent results of N. Burq,
P. G\'{e}rard and N. Tzvetkov (\cite{BGT1}). In fact, Tataru showed
that the estimate $p_\g(\phi_i)\leq C_\g\mu_i^{1/6}$ holds for any
smooth non-flat curve $\g\subset Y$ and $p^\g(\phi_i)\leq
C_\g\mu_i^\qtr$ for a flat smooth curve. For general $Y$ the
second bound is sharp. Namely, one can see that on the standard
sphere for each eigenvalue $\mu=n(n+1)$ and for each geodesic
$\ell$ (i.e. an equator) there exists an eigenfunction
$Y_{\mu,\ell}$ (e.g. the lowest associated spherical harmonic
$Y_n^n$, see \cite{Ma}) of the $L^2$-norm one such that
$\int_\ell|Y_{\mu,\ell}|^2d\ell$ is of order $\mu^\qtr$.

Moreover, further results in \cite{BGT2} greatly extended our results to $L^p$-norms and made the present paper obsolete.  We also mention \cite{Bo} and \cite{So} among further developments which reversed
\cite{BGT2} and showed the realtion to $L^4$-norms of eigenfunctions on $Y$.

2. We conjecture that the bound $p^{\ell/\s}(\phi_i)\ll\mu^\eps$
holds. This is consistent with the conjecture of P. Sarnak
\cite{Sa1} claiming that $\sup|\phi_i|\ll \mu^\eps$.
Unfortunately, for a general compact hyperbolic surface $Y$ the
only known improvement in the bound for the supremum of an
eigenfunction is logarithmic: $\sup|\phi_i|\leq
C\mu_i^{1/4}/\ln\mu_i$ (\cite{Be}).

We also conjecture that bounds
$|a_n^{\s/\ell}(\phi_i)|\ll\max\{n^\eps,\mu_i^\eps\}$ hold for any
$\eps>0$.

 For Riemann surfaces of the number-theoretic origins and
the special basis of eigenfunctions of number-theoretic
significance (the so-called Hecke-Maass basis) H. Iwaniec and P.
Sarnak \cite{ISa} have improved the exponent $1/4$ in the supremum
norm bound above to the exponent $5/24$. Corolaries A and B show
that the quantities $p^{\ell/\s}(\phi_i)$ are more accessible than
the supremum of eigenfunctions.

We also note that for 3-dimensional hyperbolic manifolds the bound
$\sup|\phi_i|\ll \mu^\eps$ does not hold in general as was shown
by Z. Rudnick and P. Sarnak \cite{RS} and hence one can not expect
that the corresponding generalized periods on geodesic spheres are
small for all eigenfunctions. It is interesting to study
corresponding periods along closed geodesics in order to see if
these are small with respect to the eigenvalue.

3. The generalized periods of the restriction of an eigenfunction
along a closed geodesic and along a geodesic circle are of utmost
interest in number theory (see \cite{Du}, \cite{KSa}). These are
the coefficients $a_\chi$ defined in \ref{method} (see
\eqref{coef-a-def} and \eqref{coef-a-def-K} for the exact
definition). Our main underlying results (Theorems A and B) give
sharp mean value bound on these coefficients. These results are
similar to the bound of G. Hardy on the average size of Fourier
coefficients of cusp forms. We formulate a conjecture concerning
the size of these periods (Conjecture \ref{CONJ}) which should be
viewed as an analog of the Ramanujan-Petersson conjecture on
Fourier coefficients of cusp forms. In the special case of
Hecke--Maass forms periods along closed geodesics and values at
some special points (i.e. Heegner points) give rise to special
values of $L$-functions (see \cite{KSa} and references therein).
In these cases using so-called convexity bounds on these
$L$-functions one obtains {\it stronger} bounds on periods along
these special curves then the bounds which trivially follow from
Theorems A and B. It is an intriguing question if there exists a
connection between  generalized periods and special values of
$L$-functions in general.

4. In this paper we treat restrictions to a {\it fixed} closed
geodesic. It is a deep problem to understand the dependence of the
norm of these restrictions and the corresponding periods on the
geodesic (e.g. to study the dependence of the constant $C_\ell$ on
the geodesic $\ell$).

The similar question for geodesic circles is of utmost interest in
the spectral theory. In particular, one would like to study the
behavior of $p_\s(\phi_i)$ for {\it small} geodesic circles with
the radius of order $\mu_i^{-\al}$ for $\haf\leq\al<0$. This is
related to so-called doubling constant of $\phi_i$ (see
\cite{NPS}).

5. Our methods work for non-compact hyperbolic surfaces of finite
volume as well. In fact, for non-compact Riemann surfaces of
finite volume one can pose a similar question for the size of
restrictions to a horocycle instead of a geodesic. Namely, let
$\eta\in Y$ be a fixed horocycle (there is a continuous family of
horocycles associated to each cusp of $Y$). Similarly to
$p^\ell(\phi_i)$, one defines $p^\eta(\phi_i)$ as the $L^2$-norm
of the corresponding restriction of the eigenfunction $\phi_i$ to
$\eta$ and the periods $p^\eta_n(\phi_i)$. Under the appropriate
normalization the coefficients $p^\eta_n(\phi_i)$ are equal to the
usual Fourier coefficients of Maass forms. One can show that for a
cusp form $\phi_i$ the bound $p^\eta(\phi_i)\leq
C_\eta\mu_i^{\frac{1}{6}}$ holds. Moreover, using nontrivial
bounds on Fourier coefficients of cusp forms for $Y$ arising from
congruence subgroups one can show that for the Hecke-Maass basis
of eigenfunctions on such surfaces the sharp bound
$p^\eta(\phi_i)\leq C_\eta\mu_i^\eps$ holds for any $\eps>0$. One
expects that the similar bound holds for general $Y$. In order to
improve the exponent $1/6$ above one have to improve over the
subconvexity bound for Fourier coefficients of Maass forms
obtained in \cite{BR1}. The analog of Theorems A and B for Fourier
coefficients of Maass forms was proved by A. Good \cite{Go}
following Hardy's method.

The paper is organized as follows. In Section \ref{reps} we remind
the notion of an automorphic representation and the correspondence
between eigenfunctions and automorphic representations. In Section
\ref{geodesics} we introduce non-negative Hermitian forms
associated with each closed geodesic and prove basic inequalities
for these forms. We also prove our main technical result, Theorem
\ref{a-bound-thm}, and formulate Conjecture \ref{a-conj} on the
size of the corresponding periods. In Section \ref{K-restrict} we
compute spectral density of a $K$-fixed vector and apply this to
prove Theorems A. In Section \ref{K-period} we prove the bound in
Theorem B for restrictions to geodesic circles along similar
lines.


{\bf Acknowledgments.} This paper is a part of a joint project
with J. Bernstein whom I would like to thank for numerous fruitful
discussions. It is a great pleasure to thank L. Polterovich for
enlightening discussions on the subject and for bringing to my
attention the paper \cite{Ta}. I also would like to thank N. Burq
for turning my attention to results in \cite{BGT1}.

The research was partially supported by BSF grant, Minerva
Foundation and by the Excellency Center ``Group Theoretic Methods
in the Study of Algebraic Varieties'' of the Israel Science
Foundation, the Emmy Noether Institute for Mathematics (the Center
of Minerva Foundation of Germany).


\section{Representation theory and eigenfunctions} \label{reps}

 It has been understood since the seminal works of A. Selberg
\cite{Se} and I. Gel'fand, S.Fomin \cite{GF} that representation
theory plays an important role in the study of eigenfunctions
$\phi_i$. Central for this role is the correspondence between
eigenfunctions of Laplacian on $Y$ and unitary irreducible
representations of the group $\PGLR$ (or what is more customary of
$\PSLR$). This correspondence allows one, quite often, to obtain
results that are more refined than similar results for the general
case of a Riemannian metric of variable curvature.

We remind the basic setting for the theory of automorphic
functions (see the excellent source \cite{GGPS} for the
representation-theoretic point of view we adopt here and \cite{Iw}
for a more classical approach based on harmonic analysis on the
upper half plane).

\subsection{Automorphic representations} \label{ureps}

We start with the geometric construction which allows one to pass
from analysis on a Riemann surface to representation theory.

One stars with with the upper half plane $\uH$ equipped with the
hyperbolic metric of constant curvature $-1$ (or equivalently one
might work with a more "homogeneous" model of the Poincar\'{e}
unit disk $D$; the use of $\uH$ is more customary in the theory of
automorphic functions). The group $\SLR$ acts on $\uH$ by the
standard fractional linear transformations. This action allows one
to identify the group $\PSLR$ with the group of all orientation
preserving motions of $\uH$. For reasons explained bellow we would
like to work with the group $G$ of all motions of $\uH$; this
group is isomorphic to $\PGLR$. Hence throughout the paper we
denote $G=\PGLR$.

Let us fix a discrete co-compact subgroup $\G \subset G$ and set
$Y=\G \sm \uH$. We consider the Laplace operator on the Riemann
surface $Y$ and denote by $\mu_i$ its eigenvalues and by $\phi_i$
the corresponding normalized eigenfunctions.

  The case when $\G$ acts freely on $\uH$  precisely corresponds to the
case discussed in \ref{M-forms} (this follows from the
uniformization theorem for the Riemann surface $Y$). Our results
hold for general co-compact subgroup $\G$ (and in fact for any
lattice $\G \subset G$).

We will identify the upper half plane $\uH$ with $G / K$, where $K
= PO(2)$ is a maximal compact subgroup of $G$ (this follows from
the fact that $G$ acts transitively on $\uH$ and the stabilizer in
$G$ of the point $z_0 = i \in \uH$ coincides with $K$).

We denote by $X$ the compact quotient $\G\sm G$ (we call it the
automorphic space). In the case when  $\G$ acts freely  on $\uH$
one can identify the space $X$ with the bundle $S(Y)$ of unit
tangent vectors to the  Riemann surface $Y = \G \sm \uH$.

 The group $G$ acts on $X$ (from the right) and hence on the space of
functions on $X$. We fix the unique $G$-invariant measure $\mu_X$
on $X$ of total mass one. Let $L^2(X)=L^2(X,d\mu_X)$ be the space
of square integrable functions and $(\Pi_X, G, L^2(X))$ the
corresponding unitary representation. We will denote by $P_X$ the
Hermitian form on $L^2(X)$ given by the scalar product. We denote
by $||\ ||_{X}$ or simply $||\ ||$ the corresponding norm and by
$\langle f,g \rangle_X$ the corresponding scalar product.

The identification  $Y=\G\sm \uH\simeq X/K$ induces the embedding
$L^2(Y)\sub L^2(X)$.
 We will always identify the space $L^2(Y)$ with the subspace of
$K$-invariant functions in  $L^2(X)$.

Let  $\phi$ be a normalized eigenfunction of the Laplace-Beltrami
operator on $Y$. Consider the closed $G$-invariant subspace
$L_\phi\sub L^2(X)$ generated by $\phi$ under the  action of $G$.
It is well-known that $(\pi,L)=(\pi_\phi, L_\phi)$ is an
irreducible unitary representation of $G$ (see \cite{GGPS}).

Usually it is more convenient to work with the space $V = L^\8$ of
smooth vectors in $L$. The unitary Hermitian form $P_X$ on $V$ is
$G$-invariant.

A smooth representation $(\pi, G, V)$ equipped with a positive
$G$-invariant Hermitian form $P$ we will call  a {\it smooth
pre-unitary representation}; this simply means that $V$ is the
space of smooth vectors in the unitary representation obtained
from $V$ by completion with respect to $P$.

   Thus starting with an automorphic function $\phi$ we constructed
   an irreducible smooth pre-unitary representation $(\pi, V)$.
   In fact we constructed this space together with a canonical
   morphism $\nu : V  \to C^\8 (X)$ since $C^\8(X)$ is the smooth
   part of $L^2(X)$.

   \defn {enhanced}A smooth pre-unitary representation
   $(\pi, G, V)$ equipped with a $G$-morphism $\nu: V \to C^\8(X)$
   we will call an {\it $X$-automorphic representation}.

   We will assume that the morphism $\nu$ is normalized,
    i.e. it carries the standard $L^2$ Hermitian
   form $P_X$ on $C^\8(X)$ into Hermitian form $P$ on $V$.

   Thus starting with an automorphic function $\phi$ we
   constructed

   (i) an $X$-automorphic irreducible pre-unitary representation
   $(\pi, V, \nu)$,

   (ii) a $K$-invariant unit vector $e_V \in V$
   (this vector is just our function $\phi$).

   Conversely, suppose we are given an irreducible smooth
   pre-unitary
    $X$-automorphic representation $(\pi,V, \nu)$ of
the group $G$ and a $K$-fixed unit vector $e_V \in V$. Then the
function $\phi = \nu(e_V) \in C^\8(X)$ is $K$-invariant and hence
can be considered as a function on $Y$. The fact that the
representation $(\pi, V)$ is irreducible implies that $\phi$ is an
automorphic function, i.e. an eigenfunction of Laplacian on $Y$.

 Thus we have established a natural correspondence between
 Maass forms $\phi$ and tuples $(\pi, V, \nu,
 e_V)$,   where $(\pi, V,\nu)$ is an $X$-automorphic irreducible
 smooth pre-unitary representation and $e_V \in V$ is a unit
 $K$-invariant vector.

It is well known that for $X$ compact the representation $(\Pi_X,
G, L^2(X))$ decomposes into a direct (infinite) sum
\begin{equation}\label{spec-L2X}
L^2(X)=\oplus_j (\pi_j, L_j)
\end{equation}
of irreducible unitary representations of $G$ (all representations
appear with finite multiplicities (see \cite{GGPS})). Let $(\pi,
L)$ be one of these irreducible "automorphic" representations and
$V = L^\8$ its smooth part. By definition $V$ is given with a
$G$-automorphic isometric morphism $\nu: V \to C^\8(X)$, i.e. $V$
is an $X$-automorphic representation.

If $V$ has a $K$-invariant vector it corresponds to a Maass form.
There are  other spaces in this decomposition which  correspond to
discrete series representations. Since they are not related to
Maass forms we will not study them in more detail.

\subsection{Representations of $\PGLR$}\label{irrep}
 All irreducible
unitary representations of $G$ are classified. For simplicity we
consider only those with a nonzero $K$-fixed vector (so called
representations of class one) since only these representations
arise from Maass forms. These are the representations of the
principal and the complementary series and the trivial
representation.

We will use the following standard explicit model for irreducible
smooth representations of $G$.

For every complex number $\lm$ consider the space $V_\lm$ of
smooth even homogeneous functions on $\br^2\sm0$ of homogeneous
degree $\lm-1 \ $ (which means that $f(ax,ay)=|a|^{\lm-1}f(x,y)$
for all $a\in\br \setminus 0$). The representation $(\pi_\lm,
V_\lm)$ is induced by the action of the group $\GLR$ given by
$\pi_\lm (g)f(x,y)= f(g\inv (x,y))|\det g|^{(\lm-1)/2}$. This
action is trivial on the center of $\GLR$ and hence defines a
representation of $G$. The representation $(\pi_\lm, V_\lm)$ is
called {\it representation of the generalized principal series}.

When $\lm=it$ is  purely imaginary  the representation $(\pi_\lm
,V_\lm)$ is pre-unitary; the $G$-invariant scalar product in
$V_\lm$ is given by $\langle f,g \rangle_{\pi_\lm}=\frac{1}{2\pi}
\int_{S^1} f\bar g d\te$. These representations are called
representations of {\it the principal series}.

 When $\lm\in (-1,1)$ the representation $(\pi_\lm ,V_\lm)$ is called
a representation of the complementary series. These
representations are also pre-unitary, but the formula for the
scalar product is more complicated (see \cite{G5}). We will not
discuss these since there are only finitely many such
representations for each $\G$ and we are interested in properties
of eigenfunctions $\phi_i$ as $\mu_i\to\8$.

 All these representations have $K$-invariant vectors.
 We fix a $K$-invariant unit vector $e_{\lm} \in V_\lm$ to be
  a function which is one on the unit circle in $\br^2$.

Representations of the principal and the complimentary series
exhaust all nontrivial irreducible pre-unitary representations of
$G$ of class one (\cite{G5}, \cite{L}).

Suppose we are given a class one $X$-automorphic representation $\nu:
V_{\lm} \to C^\8(X)$; we assume $\nu$ to be an isometric
embedding. Such $\nu$ gives rise to an eigenfunction of the
Laplacian on the Riemann surface $Y = X/K$ as before. Namely, if
$e_{\lm} \in V_\lm$ is a unit  $K$-fixed vector then the function
$\phi = \nu(e_\lm)$ is a normalized eigenfunction of the Laplacian
on the space $Y = X/K$ with the eigenvalue
$\mu=\frac{1-\lm^2}{4}$. This explains why $\lm$ is a natural
parameter to describe Maass forms.

\section{Closed geodesics, restrictions and periods}\label{geodesics}

\subsection{Closed geodesics}\label{geod}
We use the well-known description of closed geodesics in terms of
hyperbolic conjugacy classes in $\G$ (see \cite{Iw}). The set of
geodesics on $\uH$ consist of semicircles centered on the absolute
($Im(z)=0$) and of vertical lines. From this and the presentation
of $Y$ as the quotient $\G\sm\uH$ it follows that for each closed
geodesic $\ell\subset Y$ there exist an associated hyperbolic
element $\g\in\G$, defined up to the conjugacy in $\G$, such that
$\g$ stabilizes an appropriate geodesic
$\widetilde{\ell}\subset\uH$ and gives rise to the one-to-one
projection $\{\g\}\sm\widetilde{\ell}\to\ell$ under the map
$\uH\to \G\sm\uH\simeq Y$. Let $\g$ be such an element. We denote
$\G_\ell\subset \G$ the cyclic subgroup generated by $\g$. The
subgroup $\G_\ell$ is defined up to the conjugacy in $\G$. Under
the described above correspondence simple closed geodesics
correspond to conjugacy classes of {\it primitive} cyclic
hyperbolic subgroups which are generated by primitive hyperbolic
elements (i.e. those $\g'\in\G$ satisfying $\g'\not=\g^n$ for any
$\g\in \G$).

We use the well-known reformulation of the above description of
closed geodesics in $Y$ in terms of closed orbits of the group of
diagonal matrices $A\subset G$ acting on $X$ on the right. Let
$\G_\ell$ be as above and let $\g\in\G_\ell$ be its generator.
Since $\g\in G$ is a hyperbolic element there exists an element
$g_\g\in G$ such that $g_\g\inv\g g_\g=a_\g\in A$. We denote
$A_\g$ the subgroup generated by $a_\g$. The orbit
$\fO_\ell=g_\g\cdot A\subset X$ of $A$ is a closed orbit which is
homeomorphic to $A/A_\g$. Under the natural map $X\to Y$ the orbit
$\fO_\ell$ is mapped one-to-one onto the closed geodesic $\ell$ we
start with. We denote by $d\fO_\ell$ the unique $A$-invariant
measure on $\fO_\ell$ of the total mass one (a more geometric
normalization which corresponds to the length of the geodesic
$\ell$ would be $\int_{A/A_\g}1d^\times a$ ). We have hence the
relation
$p^\ell(\phi)=length(\ell)\int_{\fO_\ell}|\phi|^2d\fO_\ell$.

\subsection{Restrictions and periods}\label{restrict}
Let $\ell$ be a fixed closed geodesic and $\g\in\G_\ell$,
$\fO=\fO_\ell$ as above. Any such $A$-orbit $\fO$ gives rise to
natural Hermitian forms and a set of functionals on the space of
functions $C^\8(X)$. These will be our main tools in what follows.

We define $H_\fO$ to be the non-negative Hermitian form on
$C^\8(X)$ given by
\begin{align}\label{H-O-form}
H_\fO(f,g)=\int_\fO f(o)\bar{g}(o)d\fO
\end{align}
for any $f,g\in C^\8(X)$. We will use the shorthand notation
$H_\fO(f)=H_\fO(f,f)$.

Let $\tilde{A_\g}$ be the set of characters $\chi:A\to S^1\subset
\bc^\times$ trivial on $A_\g$. This is an infinite cyclic group
generated by a character $\chi_1$. Hence
$\tilde{A_\g}=\{\chi_n=\chi_1^n,\ n\in\bz\}$. We fix a point
$\dot{o}\in\fO$. To a character $\chi\in\tilde{A_\g}$ we associate
the function $\chi_.:\fO\to S^1$ given by
$\chi_.(\dot{o}a)=\chi(a)$.  For a character $\chi\in
\tilde{A_\g}$ we define the functional on $C^\8(X)$ given by

\begin{align}\label{O-funct}
d^{aut}_{\chi,\fO}(f)=\int_\fO f(o)\bar{\chi}_.(o)d\fO
\end{align}
for any $f\in C^\8(X)$. The functional $d^{aut}_{\chi,\fO}$ is
$\chi$-equivariant:
$d^{aut}_{\chi,\fO}(R(a)f)=\chi(a)d^{aut}_{\chi,\fO}(f)$ for any
$a\in A$, where $R$ is the right action of $G$ on the space of
functions on $X$. For a given orbit $\fO$ and a choice of a
generator $\chi_1$ we will use the shorthand notation
$d^{aut}_{n}=d^{aut}_{\chi_n,\fO}$.

The functions  $\{(\chi_n)_.\}$ form an orthonormal basis for the
space $L^2(\fO,d\fO)$. Hence we have the standard Plancherel
formula:
\begin{align}\label{Plancherel-1}
H_\fO(f)=\sum_n |d^{aut}_n(f)|^2\ .
\end{align}

Let $V$ be an automorphic representation. We consider the
non-negative Hermitian form $H^V_\fO$ on $V$ given by the
restriction of $H_\fO$ and non-negative Hermitian forms of rank
one $Q_n^{aut}(\cdot)=|d^{aut}_n(\cdot)|^2$ restricted to $V$. We
rewrite \eqref{Plancherel-1} as
\begin{align}\label{Plancherel-2}
H^V_\fO=\sum_n Q^{aut}_n
\end{align}
and view this as an equality of non-negative Hermitian forms on
$V$.

\subsubsection{Geometric inequality for $H_\fO$}
The Hermitian form $H_\fO$ is defined through the integral over a
compact set in $X$ and hence its average over the action of $G$ is
bounded by the standard Hermitian form $P_X$ on $L^2(X)$. Namely,
the group $G$ naturally acts on the space of Hermitian forms on
$C^\8(X)$. We denote this action by $\Pi$. We extend $\Pi$ to the
action of the algebra $H(G)=C^\8_c(G,\br)$ of smooth real valued
functions with compact support. We have the following basic

\begin{lem}{L-2-bound} For any $h\in H(G)$, $h\geq 0$ there
exists a constant $C=C_h$ such that $$\Pi(h)H_\fO\leq CP_X\ .$$
\end{lem}
\begin{proof} Let $u\in C^\8(X)$. Then $P_X(u)=<\mu_X,|u|^2>$ and
$\Pi(h)H_\fO(u)=<\mu',|u|^2>$, where $\mu'=\Pi(h)d\fO$. The
measure $\mu'$  is smooth (since $h$ is smooth) and $X$ is compact
hence the measure $\mu'$ is bounded by $C\mu_X$.
\end{proof}

\subsection{Homogeneous functionals}\label{homo-funct}
\subsubsection{Uniqueness of homogeneous functionals}\label{unique}
The functionals $d^{aut}_n$ introduced above are $\chi_n$-
equivariant. The central fact about such functionals is the
following uniqueness result:

\begin{thm}{ubi} Let $(\pi,V)$ be an irreducible smooth admissible
representations of $G$ and let $\chi:A\to \bc^\times$ be a
multiplicative character. Then $\dim\Hom_A(V,\chi)= 1$.
\end{thm}

The uniqueness statement is a standard fact in the representation
theory of $G$. It easily follows from the existence of the
Kirillov model (see \cite{GGPS}). There is no uniqueness of
trilinear functionals for representations of $\SLR$ (the space is
two-dimensional: it splits into even and odd functionals). This is
the reason why we prefer to work with $\PGLR$ (although our method
could be easily adopted to $\SLR$).

\subsubsection{Model homogeneous functionals} \label{modfunc}

For every $\lm \in \bc$ we denote by $(\pi_\lm,V_\lm)$  the smooth
class one representation of the generalized principle series of
the group $G=\PGLR$ described in \ref{irrep}. We will use the
realization of $(\pi_\lm, V_\lm)$ in the space of smooth
homogeneous functions on $\br^2 \setminus 0$ of homogeneous degree
$\lm - 1$ (see \cite{G5}).

For explicit computations it is often convenient to pass from
plane model to a line model. Namely, the  restriction of functions
in $V_\lm$ to the line $(x,1) \subset \br^2$ defines an
isomorphism of the space $V_\lm$ with the space $C^\8(\br)$ of
even smooth functions on $\br$ decaying on infinity as
$|x|^{\lm-1}$ so we can think about vectors in $V_\lm$ as
functions on $\br$. Under such an identification the action of the
diagonal subgroup is given by
\begin{equation}
 \pi_\lm(a)f(x,1)=f(a\inv x,
a)=|a|^{\lm-1}f(a^{-2}x,1)\ .
\end{equation}

Where we used the shorthand notation
$a=diag(a,a\inv)=\left(\begin{array}{cc}a&  \\
& a\inv\end{array}\right)$.

Note that in this model a $K$-fixed unit vector is given by
$e_0^\lm(x)=c(1+x^2)^{(\lm-1)/2}$ (where the normalization
constant $c$ is independent of $\lm$).

Let $s\in\bc$ and $\chi_s:A\to\bc^\times$, $\chi_s(a)=|a|^s$ be
the corresponding character. From the description of the action of
$A$ in the line model we see that the functionals (distributions)
$d^{mod}_{s,\lm}$ defined by the kernel $|x|^{-\haf-\lm/2+s/2}$
are $(A,\chi_s)$-equivariant functionals on $V_\lm$.
 Namely, we define
\begin{align}\label{d-s-mod}
d^{mod}_{s,\lm}(v)=\int |x|^{-\haf-\lm/2+s/2}v(x)dx\ .
\end{align}
We have
$d^{mod}_{s,\lm}(\pi_\lm(a)v)=\chi_s(a)d^{mod}_{s,\lm}(v)$. In
particular an $A$-invariant functional is given by
\begin{align}\label{invar-mod}
d^{mod}_{0,\lm}(v)=\int |x|^{-\haf-\lm/2}v(x)dx\ .
\end{align}

We note that for the general value of $s\in\bc$ one have to
understand the above integrals in a regularized sense but we will
be interested in unitary characters only ($s\in i\br$) for which
the integrals above are absolutely convergent.

\subsubsection{Coefficients of proportionality} The uniqueness of
homogeneous functionals implies that $d^{aut}$ and $d^{mod}$ are
proportional. Namely, let  $\fO$ be a closed  $A$-orbit, $A_\g$
the corresponding  subgroup with a generator
$a_\g=diag(a_\g,a_\g\inv)$ and let $V$ be an automorphic
representation isomorphic to a representation of the principal
series $V_\lm$. We denote by $q=q_\g=1/\ln(a_\g)$. For  $n\in\bz$
we consider the kernel $|x|^{-\haf-\lm/2+inq}$ and the
corresponding homogeneous functional $d^{mod}_{n,\lm}$ defined on
$V_\lm$ by means of this kernel. The set $\{d^{mod}_{n,\lm},\
n\in\bz\}$ exhaust the set of all $\chi$-invariant functionals on
$V$ as $\chi\in\tilde{A}_\g$. It follows from the uniqueness
theorem that for any $n\in\bz$ and any automorphic representation
$V$ which is isomorphic to $V_\lm$ there exists a constant
$a_{n,V}\in \bc$ such that
\begin{align}\label{coef-a-def}
d^{aut}_{\chi_n,\fO}=a_{n,V}\cdot d^{mod}_{n,\lm}\ .
\end{align}
The constant $a_{n,V}$ depends on the parameter $\lm$ and also on
the isometry $\nu_V:V\to V_\lm$ (e.g. when the multiplicity of the
corresponding eigenvalue of $\Dl$ is greater than one). We will,
however, use the notation $a_{n,\lm}=a_{n,V}$ suppressing this
difference as our method is not sensitive to the multiplicity of
$V$.

We denote  by $Q^{mod}_{n,\lm}$ the non-negative Hermitian form
$Q^{mod}_{n,\lm}(\cdot)=|d^{mod}_{n,\lm}(\cdot)|^2$. Taking into
account \eqref{Plancherel-2} we arrive at our basic relation
\begin{align}\label{Plancherel-3}
H^V_\fO=\sum_n |a_{n,\lm}|^2Q^{mod}_{n,\lm}\
\end{align}
of non-negative Hermitian forms on $V$.

\subsection{Average bound} We formulate now our main result (Theorem B):

\begin{thm}{a-bound-thm} For a given orbit $\fO$ as above there
exists a constant $C=C_\fO$ such that for any $T\geq 1$ the
following bound holds
\begin{align}\label{a-bound}
\sum_{ |n|\leq T} |a_{n,\lm}|^2\leq C\cdot \max(T,|\lm|)\ .
\end{align}
for any automorphic representation which is isomorphic to a
representation of the principal series $V_\lm$.
\end{thm}
\subsubsection{Test functionals} The proof of Theorem
\ref{a-bound-thm} is based on the notion of positive test
functionals on the space $\mathcal{H}(V)$ of Hermitian forms on an
automorphic representation $V$. Let
$\mathcal{H}(V)^+\subset\mathcal{H}(V)$ be the set of nonnegative
Hermitian forms.

\begin{defn}{def} A positive functional on the space $\mathcal{H}(V)$
is an additive map $\rho:\mathcal{H}(V)^+\to\br^+\bigcup \8$.
\end{defn}

The  basic example of such functionals is the functional
$\rho_v(H)=H(v)$ defined for any vector $v\in V$.

We construct below a special family of positive functionals
$\rho_T$ parameterized by the real parameter $T\geq 1$ in order to
bound coefficients $a_{n,\lm}$.

\begin{prop}{prop} There exist a constant $C$ such that for any $T\geq
1$ and an automorphic representation $V\simeq V_\lm$ we can find a
positive functional $\rho_T$ on $\mathcal{H}(V)$  satisfying
\begin{align}\label{prop-1}
\rho_T(H^V_\fO)\leq CT
\end{align}
\begin{align}\label{prop-2}
\rho_T(Q^{mod}_{n,\lm})\geq 1\ {\rm for\ any}\ |n|,|\lm|\leq T
\end{align}
\end{prop}

\begin{proof} We construct the positive functional $\rho_T$ by integrating
an elementary positive functional $\rho_{v_T}$ with the specially
chosen vector $v_T\in V_\lm$ against a smooth compactly supported
non-negative function $h\in H(G)$. Namely, let $\delta$ be a
smooth non-negative function with the support
$supp(\delta)\subset[-0.1,0.1]$ and satisfying
$\int\delta(x)dx=1$. For any real $T$ consider the function
$v_T(x)=T\cdot\delta(T(x-1))$ and view it as a vector in $V_\lm$.
We have
\begin{align}\label{norm-1}
P_X(v_T)=\|v_T\|^2_{V_\lm}=c_1T\ .
\end{align}
 We also have
 \begin{align}\label{big-trace}
Q^{mod}_{n,\lm}(v_T)\geq c_2
\end{align}
 for
$\max(|n|,|\lm|)\leq T$ since $v_T$ has the support on the
interval of the size smaller than $T$ around $1$ and the kernel of
$Q^{mod}_{n,\lm}$ is given by the oscillating function with the
phase having the variation on the support of $v_T$ smaller than
$\haf$.

Let $U\subset G$ be a small fixed neighborhood of the identity
such that $g\inv\cdot[0.9,1.1]\subset [\haf,\frac{3}{2}]$ for any
$g\in U$ under the standard action of $G$ on $\br$. For any $g\in
U$ the function $\pi_\lm(g)v_T$  is very similar to the original
function $v_T$ and supported off the singularities of the kernel
of $Q^{mod}_{n,\lm}$. In particular, $\pi_\lm(g)v_T$ satisfies
conditions \eqref{norm-1} and \eqref{big-trace}, possibly with
other constants. Hence the positive functional $\Pi(h)\rho_{v_T}$
satisfies the same conditions for any $h\in H(V)$ with $supp(h)\in
U$. Taking appropriate $h$ we see that $\rho_T=\Pi(h)\rho_{v_T}$
satisfies the condition \eqref{prop-2}. Moreover, taking into
account Lemma \ref{L-2-bound} we see that the condition
\eqref{prop-1} is also satisfied.
\end{proof}
\subsubsection{Proof of Theorem \ref{a-bound-thm}}
Taking the positive functional $\rho_T$ from Proposition
\ref{prop} we obtain:
\begin{align}
CT\geq\rho_T(H^V_\fO)=\sum_n|a_{n,\lm}|^2\rho_T(Q^{mod}_{n,\lm})\geq
\sum_{|n|\leq T}|a_{n,\lm}|^2\rho_T(Q^{mod}_{n,\lm})\geq \haf
\sum_{|n|\leq T}|a_{n,\lm}|^2
\end{align}
for any $T\geq |\lm|$.

\rightline{$\Box$}
\subsubsection{A conjecture} \label{a-conj}
While the bound in Theorem \ref{a-bound-thm} is sharp, for a
single coefficient in the sum it gives only $|a_{n,\lm}|\ll
\max(|n|^\haf,|\lm|^\haf)$. Such a bound follows easily from the
bound on the supremum norm  on automorphic representations
obtained in \cite{BR2}. This situation is reminiscent of a sharp
bound on the average size of the standard Fourier coefficients of
Maass forms versus bounds on a single Fourier coefficient towards
Ramanujan-Peterson conjecture. On the basis of this analogy we
propose the following

\begin{conj}{CONJ}
 For a fixed closed orbit $\fO$ the coefficients $a_{n,\lm}$
 satisfy the following bound
\begin{align}\label{coef-a-bound-conj} |a_{n,\lm}|\ll
(\max(|n|,|\lm|))^\eps
\end{align}
for any $\eps>0$.
\end{conj}

\section{Restriction to closed geodesics}\label{K-restrict}
In this section we use our main technical result Theorem
\ref{a-bound-thm} in order to prove Theorem A from the
Introduction.

\subsection{Spectral density of a $K$-fixed vector}
In order to use the Plancherel relation \eqref{Plancherel-1} and
the average bound \eqref{a-bound} we need to compute spectral
density of a $K$-fixed vector $e_0^\lm\in V_\lm$ with respect to
the model forms $Q^{mod}_{n,\lm}$. We have
\begin{align}\label{K-value-mod}
b_{n,\lm}=d^{mod}_{n,\lm}(e_0^\lm)=\int
|x|^{-\haf-\lm/2+inq}(1+x^2)^{\lm/2-\haf}dx=
\frac{\G(\frac{1-\lm+inq}{4})\G(\frac{1-\lm-inq}{4})}{\G(\frac{1-\lm}{2})}
\ .
\end{align}

Where we have used the table integral (\cite{Ma}):
\begin{align}\label{table}
\int|x|^{s}(1+x^2)^{t}dx=
B\left(\frac{s+1}{2},-t-\frac{s+1}{2}\right)=
\frac{\G(\frac{s+1}{2})\G(-t-\frac{s+1}{2})}{\G(-t)} \ .
\end{align}
From the exact expression in \eqref{K-value-mod} and the Stirling
formula for the asymptotic of the $\G$-function we see that there
are constants $c_i$ such that
\begin{equation}\label{K-spec-bound}
|b_{n,\lm}|^2\leq\Biggl\{
\begin{array}{ll}
    c_1|\lm|\inv & \hbox{for\ $|nq|\leq0.9|\lm|$;} \\
    c_2|\lm|^{-\haf}, & \hbox{for\ $0.9|\lm|\leq|nq|\leq1.1|\lm|$;} \\
    c_3e^{-0.1nq}, & \hbox{$1.1|\lm|\leq |nq|$.} \\
\end{array}
\end{equation}
{\it Remark.} The constants  $0.9$ and $1.1$ could be substituted
by $1-\s$ and $1+\s$ for any $0<\s<1$.
\subsection{Proof of Theorem A}\label{proof-A}
Let $\phi_i$ be a norm one Maass form in an automorphic
representation $V\simeq V_\lm$. We deduce from
\eqref{Plancherel-1}, \eqref{a-bound} and \eqref{K-spec-bound}
that
\begin{align}
\frac{1}{{\rm length}(\ell)}p^\ell(\phi_i)=
\int_\fO|\phi_i|^2d\fO=H^V_\fO(e_0^\lm)=\sum_n|a_{n,\lm}|^2Q^{mod}_n(e^\lm_0)=\\
=\sum_{|n|\leq1.1|\lm|} |a_{n,\lm}|^2|b_{n,\lm}|^2
+\sum_{1.1|\lm|\leq|n|} |a_{n,\lm}|^2|b_{n,\lm}|^2\leq
|\lm|^{-\haf} \sum_{1.1|\lm|\leq|n|} |a_{n,\lm}|^2+C'\leq
C|\lm|^{\haf}\
\end{align}
since \eqref{a-bound} implies that $\sum_{|n|\leq
1.1|\lm|}|a_{n,\lm}|^2\leq C|\lm|$ and the summation by parts
implies that $\sum_{1.1|\lm|\leq|n|} |a_{n,\lm}|^2e^{-0.1qn}\leq
C'$. This gives the bound \eqref{main-bound} since $|\lm|\approx
\mu^\haf$.

\rightline{$\Box$}

\section{Restriction to geodesic circles}\label{K-period}
In this section we prove Theorem B on restriction to geodesic
circles. The proof goes along same lines as the proof of Theorem A
for closed geodesics.

\subsection{Geodesic circles} We fix a maximal compact subgroup
$K\subset G$ and the identification $G/K\to\uH$, $g\mapsto g\cdot
i$. Let $y\in Y$ be a point and $\pi:\uH\to \G\sm \uH\simeq Y$ the
projection as before. Let $R_y>0$ be the injectivity radius of $Y$
at $y$. For any $r\leq R_y$ we define the geodesic circle of
radius $r$ centered at $y$ to be the set $\s(r,y)=\{y'\in Y|{\rm
d}(y',y)=r\}$. Since $\pi$ is a local isometry we have that
$\pi(\s_\uH(r,z))=\s(r,y)$ for any $z\in\uH$ such that $\pi(z)=y$
where $\s_\uH(r,z)$ is a corresponding geodesic circle in $\uH$
(all geodesic circles in $\uH$ are Euclidian circles though with
the different from $y$ center). We associate to any such circle on
$Y$ an orbit of a compact subgroup in $X$. Namely, any geodesic
circle on $\uH$ is of the form $\s_\uH(r,z)=hKg\cdot i$ with $h,\
g\in G$ such that $h\cdot i=z$ and $hg\cdot i\in \s_\uH(r,z)$
(i.e. an $h$-translation of a standard geodesic circle around
$i\in\uH$ passing through $g\cdot i$). Note, that the radius of
the circle is given by the distance ${\rm d}(i,g\cdot i)$ and
hence $g\not\in K$ for a nontrivial circle. Given the geodesic
circle $\s(r,y)\subset Y$ which gives rise to a circle
$\s_\uH(r,z)\subset \uH$ and the corresponding elements $g,\ h\in
G$ we consider the compact subgroup $K_\s=g\inv Kg$ and the orbit
$\fO_\s=hg\cdot K_\s\subset X$. Clearly we have $\pi(\fO_\s)=\s$.
We endow the orbit $\fO_\s$ with the unique $K_\s$-invariant
measure $d\fO_\s$ of the total mass one (from geometric point of
view a more natural measure would be the length of $\s$). We have
then $p^\s(\phi)=length(\s)\int_{\fO_\s}|\phi|^2d\fO_\s$.

We also note that for what follows the restriction $r<R_y$ is not
essential. From now on we assume that $\fO\subset X$ is any orbit
of $K'$. The restriction $r<R_y$ implies that the projection
$\pi(\fO)\subset Y$ is a smooth non-self intersecting curve on $Y$
and is not essential for our method.

\subsection{Hermitian forms} Let $\s$, $K'=K_\s$ and $\fO=\fO_\s$ be
as above. We define the non-negative Hermitian form on $C^\8(X)$
by
\begin{align}
H_\fO(f,g)=\int_\fO f(o)\bar g(o)d\fO
\end{align}
for any $f,\ g\in C^\8(X)$. We will use the shorthand notation
$H_\fO(f)=H_\fO(f,f)$. Examining the proof of Lemma
\ref{L-2-bound} we have, in the notations of \ref{L-2-bound},

\begin{lem}{L-2-bound-K} For any $h\in H(G)$, $h\geq 0$ there
exists a constant $C=C_h$ such that $$\Pi(h)H_\fO\leq CP_X\ .$$
\end{lem}

\rightline{$\Box$}

\subsubsection{Characters} We fix a point $\dot{o}\in\fO$. To
a character $\chi:K'\to S^1$ we associate a function
$\chi_.(\dot{o}k')=\chi(k')$, $k'\in K'$ on the orbit $\fO$ and
the corresponding functional on $C^\8(X)$ given by

\begin{align}\label{O-funct-K}
d^{aut}_{\chi,\fO}(f)=\int_\fO f(o)\bar{\chi}_.(o)d\fO
\end{align}
for any $f\in C^\8(X)$. The functional $d^{aut}_{\chi,\fO}$ is
$\chi$-equivariant:
$d^{aut}_{\chi,\fO}(R(k')f)=\chi(k')d^{aut}_{\chi,\fO}(f)$ for any
$k'\in K'$, where $R$ is the right action of $G$ on the space of
functions on $X$. For a given orbit $\fO$ and the choice of a
generator $\chi_1$ of the cyclic group $\hat{K'}$ we will use the
shorthand notation $d^{aut}_{n}=d^{aut}_{\chi_n,\fO}$, where
$\chi_n=\chi_1^n$. The functions $(\chi_n)_.$ form an orthonormal
basis for the space $L^2(\fO,d\fO)$.

Let $V$ be an irreducible automorphic representation. We introduce
non-negative Hermitian forms of rank one
$Q_n^{aut}(\cdot)=|d^{aut}_n(\cdot)|^2$ restricted to $V$ and
consider the Plancherel formula restricted to $V$:
\begin{align}\label{Plancherel-2-K}
H^V_\fO=\sum_n Q^{aut}_n\ .
\end{align}

Let $V\simeq V_\lm$ be a representation of the principal series.
We have $\dim\Hom_{K'}(V_\lm,\chi)\leq 1$ for any character $\chi$
of $K'$ (i.e. the space of $K$-types is at most one dimensional).
This is well-known in representation theory of $\PGLR$ and could
be seen from the isomorphism $V_\lm\simeq C^\8_{even/odd}(S^1)$,
for example. In fact $\dim\Hom_K(V_\lm,\chi_n)=1$ iff $n$ is even
for $V$ corresponding to a Maass form.

Consider the model $V_\lm\simeq C^\8(S^1)$ and the standard
vectors (exponents) $e_n=\exp(2\pi in)\in C^\8(S^1)$ which form
the basis of $K$-types for the {\it standard} maximal compact
subgroup $K$. For any $n$ such that $\dim\Hom_K(V_\lm,\chi_n)=1$
the function $e'_n=\pi_\lm(g\inv)e_n$ defines the model functional
on $V_\lm$ through
$d^{mod}_{n,\lm}(v)=d^{mod}_{\chi_n,\lm}(v)=<v,e'_n>$ which is
$\chi_n$-equivariant with respect to $K'$. Introducing the
Hermitian forms $Q_n^{mod}=|d_{n,\lm}(\cdot)|^2$ we arrive at
basic relations
\begin{align}\label{coef-a-def-K}
Q^{aut}_{n,\lm}= |a_{n,\lm}|^2Q^{mod}_{n,\lm}\ .
\end{align}
\begin{align}\label{Plancherel-3-K}
H^V_\fO=\sum_n |a_{n,\lm}|^2Q^{mod}_{n,\lm}\ .
\end{align}

\subsubsection{Average bound} By examining the proof of Theorem
\ref{a-bound-thm} we arrive at a similar statement.

\begin{thm}{a-bound-thm-K} For a given orbit $\fO$ as above there
exists a constant $C=C_\fO$ such that for any $T\geq 1$ the
following bound holds
\begin{align}\label{a-bound-K}
\sum_{ |n|\leq T} |a_{n,\lm}|^2\leq C\cdot \max(T,|\lm|)\ .
\end{align}
for any automorphic representation which is isomorphic to a
representation of the principal series $V_\lm$.
\end{thm}

The proof again is based on the existence of appropriate test
functionals on the space of Hermitian forms on $V_\lm$ satisfying
the same conditions as in Proposition \ref{prop}.

\rightline{$\Box$}

\subsection{Spectral density} We are left to compute spectral
decomposition of the $K$-fixed vector $e_0\in V_\lm$ with respect
to the basis of $K'$-types. Namely, we need to estimate
coefficients

\begin{align}\label{K-vaues-K}
c_{n,\lm}=d^{mod}_{n,\lm}(e_0)=<e_0,e'_n>=
<e_0,\pi_\lm(g\inv)e_n>=<\pi_\lm(g)e_0,e_n>\
.
\end{align}
These are finite $K$-types matrix coefficients of the spherical
vector $e_0\in V_\lm$. These matrix coefficients satisfy the
following crude bound ( compare \eqref{K-spec-bound}).

\begin{lem}{K-spec-bound-lem} Let $g\in G$ be fixed, $g\neq e$.
There exists a constant $c=c_g> 1$ and constants $c_1,\ c_2>0$
such that
\begin{equation}\label{K-spec-bound-K}
|c_{n,\lm}|^2\leq\Biggl\{
\begin{array}{ll}
    c_1|\lm|\inv & \hbox{for\ $|2\pi n|\leq0.9c|\lm|$;} \\
    c_2|\lm|^{-\frac{2}{3}}, & \hbox{for\ $0.9c|\lm|\leq|2\pi n|\leq1.1c|\lm|$;} \\
     o(n^{-N}) & for\ any\ N>0, \hbox{for\ $1.1c|\lm|\leq |2\pi n|$.} \\
\end{array}
\end{equation}
\end{lem}
\begin{proof} We want to estimate quantities $<\pi_\lm(g)e_0,e_n>$.
Namely, we need to estimate the coefficients of the Fourier
expansion of $f_\lm(\theta)=\pi_\lm(g)e_0(\theta)\in C^\8(S^1)$.
The function $f_\lm$ is given by the formula
$f_\lm(\theta)=|g'(\theta)|^{-\frac{1+\lm}{2}}$. We consider the
corresponding oscillatory integral
\begin{align}\label{c's}
c_{n,\lm}=<f_\lm,e_n>=\int_{S^1}|g'(\theta)|^{-\haf}
e^{\haf\lm\ln|g'(\theta)|-2\pi in\theta}d\theta\ .
\end{align}
Clearly, $g'(\theta)$ is a smooth real valued non-zero function.
Let $c=\max_{S^1} d\log|g'|$. It is easy to see (for example by
taking the Cartan decomposition $g=k_1diag(a,a\inv)k_2$, where
$k_1,\ k_2\in K$) that the phase of the oscillatory integral
\eqref{c's} has non-degenerate critical points for $|2\pi n|\leq
0.9 c|\lm|$ no critical points for $|2\pi n|\geq 1.1 c|\lm|$ and a
critical point with the degeneration of a degree at most 3 (i.e.
of a type $ct^3$ in a local parameter $t$)  for $ 0.9
c|\lm|\leq|2\pi n|\leq 1.1 c|\lm|$. The bound
\eqref{K-spec-bound-K} follows now from the stationary phase
method.
\end{proof}
\subsubsection{Proof of Theorem B} Examining the proof of Theorem A
in Section \ref{proof-A}  we arrive immediately at the bound
\eqref{main-bound-K} in Theorem B.

\rightline{$\Box$}

\begin{rem}{b-ver-c}
 The discrepancy in the behavior of coefficients
$b_{n,\lm}$ and $c_{n,\lm}$ (compare \eqref{K-spec-bound} to
\eqref{K-spec-bound-K}) is a result of the difference in the type
of degeneration of the phase of the oscillatory integrals
\eqref{K-value-mod} and \eqref{c's}. This is reflected in the
difference of corresponding exponents in bounds in Corollaries A
and B. While the phase in the integral in \eqref{c's} for the
coefficients related to geodesic circles has degeneration of a
degree at most 3, the corresponding phase in \eqref{K-value-mod}
for the coefficients related to closed geodesics has the critical
point which coincides with the singularity of the amplitude. We
note that underlying average bounds on coefficients $a_{\chi,\lm}$
relating the geometric Hermitian form $H_\fO$ to model Hermitian
forms $Q^{mod}_{\chi,\lm}$ are the same in both cases (compare
Theorems \ref{a-bound-thm} and \ref{a-bound-thm-K}).
\end{rem}

\end{document}